\documentclass[11pt, a4paper]{article}
\usepackage[subpreambles=false]{standalone}%
\usepackage[english]{babel}
\usepackage[top=1.5cm, bottom=1.5cm, left=1.5cm, right=1.5cm, headheight=.5cm, headsep=.5cm, a4paper]{geometry}
\usepackage{caption}
\usepackage{subcaption}
\usepackage{xcolor}
\usepackage{authblk}
\usepackage{graphicx}
\usepackage{graphics}
\usepackage{pgfplots}
\pgfplotsset{compat=1.15}
\usepackage{mathrsfs}
\usetikzlibrary{arrows}
\usepackage[unicode, colorlinks = true]{hyperref}
\hypersetup{
colorlinks=true,
linkcolor=[rgb]{.75,.15,.10},
citecolor=[rgb]{.35,.65,.55},
urlcolor=[rgb]{.88,.33,.88}}
\usepackage{caption,subcaption,xcolor,amssymb,amsmath,amsthm}
\usepackage[shortlabels]{enumitem}
\usepackage{pgfplots}
\usepackage[normalem]{ulem}
\usepackage{soul}
\pgfplotsset{compat=1.15}
\usepackage{mathrsfs}
\usetikzlibrary{arrows}
\usepackage{dsfont}
\def\abs#1{\left\lvert#1\right\rvert}
\def\R{\mathds{R}}
\def\wt{\widetilde}
\def\wh{\widehat}
\def\ol{\overline}

\def\phi{\varphi}
\def\d{{\mathrm d}}
\newtheorem{proposition}{Proposition}

\begin{document}
\sloppy

\title{Separation of the initial conditions in the inverse problem for 1D non-linear tsunami wave run-up theory\footnote{To appear in Studies in Applied Mathematics}} 

\author[1]{Alexei Rybkin}
\author[1]{Oleksandr Bobrovnikov}
\author[2]{Noah Palmer}
\author[3]{Daniel Abramowicz}
\author[4,5]{Efim Pelinovsky}

\affil[1]{University of Alaska Fairbanks, Fairbanks, Alaska, United States}
\affil[2]{University of Colorado Boulder, Boulder, Colorado, United States}
\affil[3]{University of San Francisco, San Francisco, California, United States}
\affil[4]{HSE University, Nizhny Novgorod, Russia}
\affil[5]{Institute of Applied Physics, Nizhny Novgorod, Russia}
\maketitle

\begin{abstract}
We investigate the inverse tsunami wave problem within the framework of the 1D non-linear shallow water equations (SWE). Specifically, we show that the initial displacement $\eta_0(x)$ and velocity $u_0(x)$ of the wave can be recovered, given the known motion of the shoreline $R(t)$ (the wet/dry free boundary), in terms of the Abel transform. We demonstrate that for power-shaped inclined bathymetries, this problem admits a complete solution for any $\eta_0$ and $u_0$, provided the wave does not break.

It is important to note that, in contrast to the direct problem (also known as the tsunami wave run-up problem), where $R(t)$ can be computed exactly only for $u_0(x)=0$, our algorithm can recover $\eta_0$ and $u_0$ exactly for any non-zero $u_0$. This highlights an interesting asymmetry between the direct and inverse problems. Our results extend the work presented in \cite{Rybkin23,Rybkin24}, where the inverse problem was solved for $u_0(x)=0$. As in previous work, our approach utilises the Carrier-Greenspan transformation, which linearises the SWE for inclined bathymetries. Extensive numerical experiments confirm the efficiency of our algorithms.
\end{abstract}


\section{Introduction}
Tsunami wave run-up is often analysed within the framework of the shallow water equations (SWE), a system of non-linear partial differential equations. In these studies, it is typically assumed that the initial water displacement and velocity of the incoming wave are known, and the objective is to compute the motion of the shoreline (also known as the dry/wet free boundary). This problem is commonly referred to as the tsunami wave run-up problem.

In this work, we are concerned with the inverse problem: determining the initial conditions (i.e., the initial displacement and velocity of the wave) based on the known motion of the shoreline. This problem has long been of practical importance for tsunami forecasting and mitigation efforts (see \cite{Miyabe, Piatanesi, Pires, Voronin, Satake21} and the literature cited therein). Since tsunami waves are governed by non-linear PDEs, the inverse tsunami wave problem poses a significant challenge. Thus, it is critical to develop reasonable models where this problem can be solved effectively.

In \cite{Rybkin23}, it was demonstrated that for a sloping plane beach, the initial displacement of the wave can be recovered under the assumption of zero initial velocity. More recently, in \cite{Rybkin24}, this result was extended to bays with inclined power-shaped cross-sections. Such three-dimensional bathymetries can be parametrised by a single function of one variable (see \cite{Rybkin21}, for precise definitions and statements). However, the assumption of zero initial velocity has remained unaddressed, limiting applicability to real-life scenarios where initial velocity plays a major role (e.g., tsunamis caused by underwater landslides). In such cases, the initial displacement is often negligible instead. In \cite{Rybkin24}, this case was left as an open problem.

In this paper, we provide a complete solution to the inverse problem for inclined power-shaped bays where both the initial displacement and initial velocity are non-zero. While our approach follows the framework established in \cite{Rybkin24}, two new insights are critical to our results, which we summarise below at a non-technical level.

The first insight concerns a key difference between the direct and inverse problems. For context, as in \cite{Rybkin24}, we use the Carrier-Greenspan hodograph transformation to linearise the shallow water equations into a specific form of the linear wave equation. This allows us to apply standard techniques from mathematical physics. Relevant works in this context include those by \cite{Shimozono}, \cite{Pedersen21}, \cite{Synolakis06}, and \cite{Aydin20} on SWE in power-shaped bays and plane beaches. For numerical treatment of channels of arbitrary cross-sections, we refer the reader to works of \cite{hernandez2011shallow}, \cite{HADIARTI202374}, \cite{Welahettige}, and \cite{WANG2019124587}. For a different kind of inverse problems for shallow water equations we suggest the recent work of \cite{HAKL2025134496}. 

However, while the Carrier-Greenspan linearisation simplifies the problem, it introduces a challenge: the initial conditions in the hodograph plane are specified on a curve that depends on the initial conditions themselves. This curve becomes a straight line only if the initial velocity is zero, which is a well-documented issue in hodograph transformations (e.g., \cite{Johnson97}). In prior studies, this issue was addressed under additional assumptions (see \cite{Carrier03}; \cite{Kanoglu06}). In \cite{Nicolsky18}, we introduced the data projection method, which offers a robust solution for non-zero initial velocity. This method, further developed in \cite{Rybkin21}, approximates initial conditions in the hodograph plane using Taylor's formula along a straight line, ensuring the solutions to the original and modified problems are as close as desired (though not exact).

For the inverse problem considered here, such complications are absent. The initial conditions can be recovered exactly as long as the Carrier-Greenspan transformation remains invertible. This leads to a surprising asymmetry between the direct and inverse problems, which appears to be a new phenomenon.

Our second novel insight is both surprising and counter-intuitive. A detailed analysis of our formulas reveals that, in the hodograph plane, the new dependent variables (interpreted physically as pressure and velocity) uniquely contribute to the solution at the shoreline, a fixed point in the hodograph plane. This allows the recovery of both the initial displacement and the initial velocity from the shoreline motion, which is described by a single function, referred to as the shoreline data or run-up function. It is particularly striking that two independent initial conditions can be deduced from a single function -- a phenomenon we term the “separation of initial conditions at the shoreline”. The run-up function separates contributions from the initial displacement and initial velocity. In practical applications, when the run-up is derived from measurements, this separation can only be approximate.

The remainder of this work is organised as follows: in Section \ref{sec:model} the mathematical model is introduced; in Section \ref{sec:solToFp} we derive the solution to the direct problem. We present and discuss the phenomenon of the separation of the initial conditions at the shoreline in Section \ref{sec:sep}.  Section \ref{sec:ip} provides a complete solution to the inverse problem for inclined power-shaped bays
when both initial displacement and initial velocity are non-zero. Some remarks left outside of the main body are given in Section \ref{sec:remarks}. Finally, we provide a numerical experiment which validates our findings and concluding remarks in Sections \ref{sec:num} and \ref{sec:conc} respectively.

\section{Mathematical model}\label{sec:model}
\begin{figure}[!b]
    \begin{subfigure}{.45\textwidth}
        \input{bath_xOz.tex}
      \caption{$(x, z)$ cross-sectional view;}
    \end{subfigure}\hfill
    \begin{subfigure}{.45\textwidth}
     \begin{center}
      \includegraphics[width = \linewidth]{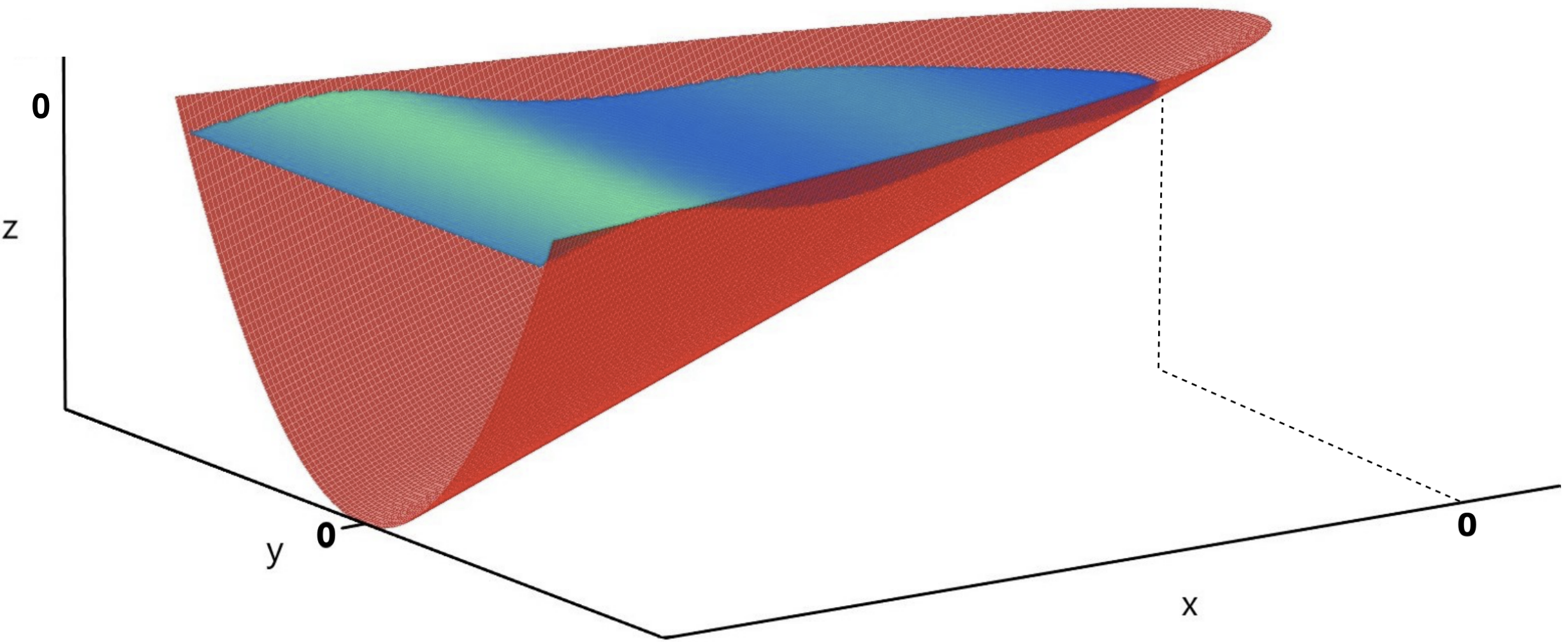}
     \end{center}
      \caption{3D view;}\label{3dview}
    \end{subfigure}
    \caption{A sketch of the parabolic bay geometry. The bathymetry \(z (x, y) = - x + \abs{y}^m \) is in red, the unperturbed water level is in yellow, an the water level is in blue. The total perturbed water depth is given by $H(x,t)=h(x)+\eta(x,t)$.}\label{fig:bath}
\end{figure}

The 1D shallow water (Saint-Venant) equations, given in dimensionless units as
\begin{equation}\label{SWE}
    \begin{aligned}
      &\partial_t\eta+u \partial_x (x+\eta) +\frac{m}{m+1}(x+\eta)\partial_x u=0,&\text{(mass)}  \\
      &\partial_t u+u\partial_x u+\partial_x \eta= 0, &\text{(momentum)}
    \end{aligned}
\end{equation}
govern water motion in an infinite parabolic-shaped bathymetry (see Figure \ref{fig:bath}) under the assumptions of no vorticity, no friction, no dispersion, and no wave breaking.
Here $u(x, t)$ is the depth averaged flow velocity over the corresponding cross-section, $\eta(x, t)$ is the water displacement exceeding the unperturbed water level, and \(m\) is the bay parameter.
The substitution
\begin{equation}\label{eq:dim}
\wt{x}= \frac{H_0 }{\alpha}{x},\quad \wt{t}=\sqrt{ \frac{H_0 }{g}}\,\frac{t}{\alpha},\quad \wt{\eta}=H_0{\eta},\quad \wt{u}=\sqrt{H_0 g}\,{u},
\end{equation}
where $H_0$ is the characteristic height of the wave, $\alpha$ is the slope of the bathymetry, and $g$ is the acceleration of gravity, turns the dimensionless system into one with dimension (dimensional variables are the ones with tildes).
There are two fundamental difficulties associated with the system \eqref{SWE}: non-linearity and the free boundary.
As a result of water motion, the shoreline is shifting; let us denote the law of \emph{vertical} motion of the shoreline as \(R(t)\).   
Typically the system \eqref{SWE} is considered together with the initial conditions
\begin{equation}\label{phys_ic}
    u(x,0) = u_0 (x) ,\quad \eta(x,0) = \eta_0(x)
\end{equation}
and the goal is to compute the vertical run-up of the wave \(R(t)\). We call this the direct or forward problem.

\section{Solution to the forward problem}\label{sec:solToFp}
In this section we outline a solution to the forward problem. The forward problem is well-studied both numerically \cite{Hartle21,Bueler22} and analytically \cite{Pelin-Maz,Kanoglu06,Rybkin14,Didenkulova11,Carrier58}. Our solution to the forward problem here mostly follows the solution of \cite{Rybkin21}.

The Carrier-Greenspan hodograph transformation (CGT), originally introduced for the sloping plane beach by \cite{Carrier58} and later generalised for power-shaped bays by \cite{Rybkin14}, allows one to linearise \eqref{SWE} and fix the free boundary. Specifically, this transform introduces new variables \(\sigma, \tau\) and new functions \(\psi,\phi\) as
\begin{equation}\label{CG}
\begin{aligned}
  &\underset{\text{\normalsize\vphantom{\LARGE X}(total  height)}}{\sigma = x + \eta(x,t)}
  ,
  &&\underset{\text{\normalsize\vphantom{\LARGE X}(delayed time)}}{\tau = t-u(x,t),}
    &\underset{\text{\normalsize\vphantom{\LARGE X}(velocity)}}{\varphi(\sigma,\tau) = u(x,t)},
    &&
    \underset{\text{\normalsize\vphantom{\LARGE X}(pressure)}}{\psi(\sigma,\tau) = \eta(x,t) + u^2(x,t)/2}.
\end{aligned}
\end{equation}
Physically, \(\sigma\) is interpreted as the total height of the wave, \(\tau\) as the delayed time, \(\phi\) as the velocity and \(\psi\) as the pressure (from Bernoulli's equation).

Like most of hodograph transformations, the CGT complicates the initial conditions. Indeed, when we set $t=0$ in the equation $\tau = t-u(x,t)$, we see that \(\tau\) is not (identically) zero unless \(u(x,0) = 0\). In the hodograph plane \((\sigma,\tau)\) the initial conditions are specified on a parametric curve which we denote \(\tau = \gamma(\sigma)\).
The system \eqref{SWE} together with \eqref{phys_ic} under \eqref{CG} becomes
\begin{equation}\label{eq:hyp}
\begin{aligned}
     & \partial_\tau \psi + \frac{m }{m+1} \sigma \partial_\sigma \phi + \phi = 0, \\
     & \partial_\tau \phi + \partial_\sigma \psi = 0,\\ 
     & \psi(\sigma, \gamma(\sigma)) = \eta_0 + u_0^2 / 2=: \psi_{\text{phys}}(\sigma), \\ 
     & \phi(\sigma, \gamma(\sigma)) = u_0 =: \phi_{\text{phys}}(\sigma).
\end{aligned}
\end{equation}
Note that this is a linear hyperbolic system. The method of data projections \cite{Rybkin21} allows one to replace the initial conditions on \(\gamma\) with initial conditions on the line \(\tau = 0\) such that the solution to the original and the solution to the new IVP are arbitrarily close. To give the reader an understanding of the complexity of the method we give the formula for the \(n\)-th order projection:
\begin{equation}
\label{eq:proj}
  \begin{pmatrix}
     \phi_n(\sigma)\\ 
     \psi_n(\sigma)
  \end{pmatrix}
  =
 \begin{pmatrix}
     \phi_{\text{phys}} (\sigma) \\ 
     \psi_{\text{phys}} (\sigma)
 \end{pmatrix}
 + 
 \sum_{k=1}^n \frac{\phi_{\text{phys}}^k(\sigma)}{k!}
 \left[
D \Delta
 \right]^k
 \begin{pmatrix}
     \phi_{\text{phys}}(\sigma)\\ 
     \psi_{\text{phys}}(\sigma)
 \end{pmatrix},
\end{equation}
where
\begin{equation}
     D =
     \begin{pmatrix}
      1 & \phi_{\text{phys}}'(\sigma) \\ 
      \frac{m \sigma}{m+1} \phi_{\text{phys}}'(\sigma) & 1
  \end{pmatrix},~
     \Delta = - \begin{pmatrix}
        0 & 1 \\ 
        \frac{m \sigma}{m+1} & 0
     \end{pmatrix}
     \partial_\sigma - 
     \begin{pmatrix}
        0 & 0 \\
        1 & 0
     \end{pmatrix}.
\end{equation}
Using \eqref{eq:proj} the IVP \eqref{eq:hyp} becomes
\begin{equation}\label{eq:hyp_proj}
    \begin{aligned}
         & \partial_\tau \psi + \frac{m }{m+1} \sigma \partial_\sigma \phi + \phi = 0, \\
         & \partial_\tau \phi + \partial_\sigma \psi = 0,\\ 
         & \psi(\sigma, 0) = \psi_{\text{proj}} (\sigma), \\ 
         & \phi(\sigma, 0) = \phi_{\text{proj}}(\sigma).
    \end{aligned}
\end{equation}
Here \(\phi_{\text{proj}} , \psi_{\text{proj}} \) are found from \(u_0 , \eta_0 \) using the infinite order data projection so that the solutions to \eqref{eq:hyp} and \eqref{eq:hyp_proj} are \emph{exactly} the same.
Note that it is impossible to find these initial conditions from the physical initial conditions exactly when solving the forward problem since Taylor's series is truncated when doing data projections.
The linear hyperbolic system \eqref{eq:hyp_proj} then can be written as a wave equation
\begin{equation} \label{eq:wave}
    \partial^2_\tau\psi = \frac{m}{m+1}\sigma\partial^2_\sigma\psi + \partial_\sigma\psi,
\end{equation}
which can be solved using the Fourier-Bessel transform. Returning back to the hyperbolic system \eqref{eq:hyp_proj} one obtains the solution
\begin{multline}\label{psisol}
    \psi(\sigma,\tau) =2 \sigma^{-\frac{1}{2m}}\int_0^\infty k
    \left\{\int_0^\infty\psi_{\text{proj}}(s)s^{\frac{1}{2m}}J_{\frac{1}{m}}(2k\sqrt{s})\,\d s\cos(\omega k\tau)
    \right. \\ \left.
    -\omega\int_0^\infty \phi_{\text{proj}}(s)s^{\frac{1}{2m}+\frac{1}{2}}J_{\frac{1}{m}+1}(2k\sqrt{s})\,\d s\sin(\omega k\tau)\right\}
    J_{\frac{1}{m}}(2k\sqrt{\sigma})\,\d k,
\end{multline}
\begin{multline}\label{phisol}
    \phi(\sigma, \tau) = \frac{2}{\omega} \sigma^{- \frac{1}{2m}-\frac{1}{2}}\int_0^\infty k
    \left\{
         \int_0^\infty \psi_{\text{proj}}(s) s^{\frac{1}{2m}} J_{\frac{1}{m}}(2 k \sqrt s)\,\d s \sin(\omega k \tau) \right.
        \\
        \left.
            + \omega  \int_0^\infty \phi_{\text{proj}} (s) s^{\frac{1}{2m}+\frac{1}{2}} J_{\frac{1}{m}+1}(2 k \sqrt s)\,\d s \cos(\omega k \tau)
    \right\}
    J_{\frac{1}{m }+1}(2 k \sqrt \sigma)\,\d k,
\end{multline}
where \(\omega = \sqrt{m / (m+1)}\) and \(J_ \nu\) is the Bessel function of the first kind of order \(\nu\).
Taking \(\sigma \to 0\) in \eqref{psisol} one obtains the relation that we call the shoreline equation
\begin{multline}\label{psistep2}
    \psi(0, \tau) =\frac{4}{\Gamma\left(1+\frac{1}{m}\right)}\int_0^\infty k^{1+\frac{1}{m}}
    \left\{\int_0^\infty\wh\psi_{\text{proj}}(\lambda)\lambda^{\frac{1}{m}+1}J_{\frac{1}{m}}(2k\lambda)\,\d\lambda\cos(\omega k\tau)
    \right.\\\left.
    -\omega \int_0^\infty \wh\phi_{\text{proj}}(\lambda)\lambda^{\frac{1}{m}+2}J_{\frac{1}{m}+1}(2k\lambda)\,\d\lambda\sin(\omega k\tau)\right\}\,\d k,
\end{multline}
where $\lambda=\sqrt{s}$, $\wh\psi_{\text{proj}}(\lambda)=\psi_{\text{proj}}(\lambda^2)$ and $\wh\phi_{\text{proj}}(\lambda)=\phi_{\text{proj}}(\lambda^2)$.
The run-up \(R(t) \) can be found from \(\psi(0,\tau)\) using \eqref{CG}, which at the shore reads as
\begin{equation}\label{CGarshore1}
    {\tau = t+\partial_t R (t),}
    \quad
    {\varphi(0,\tau) = -\partial_t R (t)}
    ,\quad
    {\psi(0,\tau) = R (t) + \tfrac{1}{2}(\partial_t R (t))^2}.
\end{equation}
Note that here we take advantage of dealing with dimensionless units: at the free boundary we have  \(\eta(x, t) = R(t)\) and \(u(x, t) = -\partial_t R(t)\) since the slope of the bathymetry is \( \pi / 4\). 
This concludes the solution to the forward problem.
\section{Separation of the initial conditions at the shoreline}\label{sec:sep}
In this section we present two propositions which describe the relationship between the pressure and velocity at the shoreline and allow for the separation of the initial conditions at the shoreline. First, like in the derivation of \eqref{psistep2}, taking \(\sigma \to 0\) in \eqref{phisol} yields
\begin{multline}\label{phistep2}
    \phi(0, \tau) = \frac{4}{\omega \Gamma(2+\frac{1}{m})}
    \int_0^\infty k^{2 + \frac{1}{m}}
    \\ 
    \left\{
          \int_0^\infty \wh \psi_{\text{proj}}(\lambda) \lambda^{1+\frac{1}{m}} J_{\frac{1}{m}}(2 k \lambda)\,\d \lambda \sin(\omega k \tau) \right.
        \\
        \left.
            +\omega \int_0^\infty \wh\phi_{\text{proj}} (\lambda) \lambda^{\frac{1}{m}+2} J_{\frac{1}{m}+1}(2 k \lambda)\,\d \lambda \cos(\omega k \tau)
    \right\}
    \,\d k.
\end{multline}
Comparing \eqref{psistep2} and \eqref{phistep2} we observe that the velocity $\phi$ at the shoreline can be re-expressed as the derivative of $\psi$. 
This leads us to the following important observation.
\begin{proposition}
    The pressure  $\psi$ and velocity $\phi$ from \eqref{CG} are related at the shoreline ($\sigma=0$) by
    \begin{equation}
\label{derivativeatshore}
\phi(0,\tau)=- \partial_\tau\psi(0,\tau).
\end{equation}
\end{proposition}
Next, we examine \eqref{psistep2}. Note that under the outer integral there is a sum of even and odd functions with respect to \(\tau\). 
This suggests that if \(\psi(0,\tau)\) could be split into its even and odd components, it would be possible to uncouple \(\wh\psi_{\text{proj}}\) and \(\wh\phi_{\text{proj}}\):
\begin{equation}\label{eq:separ}
    \begin{aligned}
        &[\psi(0, \tau)]_{\text{even}} =\frac{4}{\Gamma\left(1+\frac{1}{m}\right)}\int_0^\infty k^{1+\frac{1}{m}}
        \int_0^\infty\wh\psi_{\text{proj}}(\lambda)\lambda^{\frac{1}{m}+1}J_{\frac{1}{m}}(2k\lambda)\,\d\lambda\cos(\omega k\tau)
       \,\d k,\\
      &[\psi(0, \tau)]_{\text{odd}} =\frac{-4 \omega}{\Gamma\left(1+\frac{1}{m}\right)}\int_0^\infty k^{1+\frac{1}{m}}
       \int_0^\infty \wh\phi_{\text{proj}}(\lambda)\lambda^{\frac{1}{m}+2}J_{\frac{1}{m}+1}(2k\lambda)\,\d\lambda\sin(\omega k\tau)\,\d k.
    \end{aligned}
\end{equation}
Here we mean even/odd in the usual sense: the function \(f\) is even if \(f(x)=f(-x)\) and is odd if \(f(x)=-f(-x)\). To allow such separation we assume \(R (t)\) to be a real analytic function (that is, to have a convergent Taylor series in the entire complex plane). The physical domain of \(R(t)\) is narrower. Analyticity of \(R (t)\) together with \eqref{CGarshore1} easily implies analyticity of \(\psi(0,\tau)\) and \(\phi(0,\tau)\),  as long as the CGT is invertible.
We arrive at our second statement.
\begin{proposition}
    \label{disentangled}
    Let $R(t)$ be an analytic function
    that fits the
    shoreline run-up data. Assume that the CGT is invertible, and so the wave does not break. Then the even component of the
    pressure at the shoreline
    $\psi(0,\tau)=R(t)+(\partial_tR(t))^2/2$, where \(t =  \tau + \phi(0, \tau) = \tau - \partial_t R (t)\),  depends only on the initial pressure of the wave and the odd component depends only on the initial velocity of the wave.
\end{proposition}
It is worth mentioning that analyticity is a sufficient condition, but not necessary. To allow separation we only need a reasonable way to split a function into its even and odd components.

\section{Inverse problem}\label{sec:ip}
We now consider the full inverse problem. Given the shoreline oscillations $R(t)$, we recover the initial velocity $u_0$ and displacement $\eta_0$. Since the conversion of \(R (t)\) to \(\psi(0,\tau)\) is trivial, there are two steps in the solution to our inverse problem:
\begin{enumerate}[(1)]
    \item from \(\psi(0,\tau)\) recover \(\psi_{\text{proj}} , \phi_{\text{proj}} \);
    \item from \(\psi_{\text{proj}}\) and \(\phi_{\text{proj}} \) find \(\eta_0\) and \(u_0 \).        
\end{enumerate}

The equations in \eqref{eq:separ} can be inverted using inverse Fourier transforms together with the integral representation of the Bessel function. However, for \(\wh\phi_{\text{proj}}\) it is more convenient to separate \(\phi(0,\tau)\) in a similar manner to that done for \(\psi(0,\tau)\) to obtain an expression for \(\wh\phi_{\text{proj}}\) in terms of an integral of the even component of $\phi(0,\tau)$. From this one obtains
\begin{align}
    \wh\psi_{\text{proj}}(\lambda) &= 
    \frac{2\sqrt\pi \Gamma(1+1/m)}{\lambda \Gamma(1/m+1/2)}\int_0^{\lambda/\pi}\left( 1 - \left( \frac{\pi\xi}{\lambda} \right)^2\right)^{1/m-1/2} \Psi(\xi)~\d\xi,     \label{finalPsi}\\
    \wh \phi_{\text{proj}}(\lambda) &= \frac{2\sqrt \pi \Gamma(2+{1}/{m})}{\lambda^{2+{2}/{m}}\Gamma({3}/{2}+{1}/{m})}\int_0^{{\lambda}/{\pi}} (\lambda^2 - \pi^2 \xi^2)^{1 / m + 1 / 2}\Phi(\xi) \,\d \xi,\label{finalPhi}
\end{align}
where $\Psi(\xi):=\psi_{\text{even}}(0,q\xi)$, $\Phi(\xi):=\phi_{\text{even}}(0,q\xi)$, $q = {2\pi}/{\omega}$, and $\xi = \tau / q$. This completes the first step in the solution of the inverse problem: from the energy at the shoreline \(\psi(0,\tau)\) we recover the initial conditions \(\psi_{\text{proj}}\) and \(\phi_{\text{proj}}\) in the hodograph \((\sigma,\tau)\) plane using \eqref{finalPsi} and \eqref{finalPhi} respectively. This recovery is exact in the following sense: precisely these initial conditions on \(\tau = 0\) in the hodograph plane produce the observed \(R (t)\). In constrast, when solving the direct problem it is impossible to obtain these conditions (on \(\tau = 0\)) from \(\eta_0(x)\) and \(u_0(x) \) exactly since the  method of data projections truncates the Taylor series.

\subsection*{Exact recovery of the initial conditions in the physical plane}
Now we recover the initial conditions \(u_0\) and \(\eta_0\) in physical space. From (\ref{psisol},~\ref{phisol}) we find \(\psi(\sigma, \tau)\) and \(\phi(\sigma,\tau)\) (their values for \emph{all} values \((\sigma,\tau)\)) and then perform the inverse CGT. We first need to find the curve \(\gamma\). Note that on \(\gamma\) (equivalently at \(t=0\)) we have from \eqref{CG} that
${\gamma(\sigma)=\tau= -u(x,0) = - \phi(\sigma, \gamma(\sigma))}$. Hence, the curve $\gamma(\sigma)$ can be found by solving
\(
\gamma(\sigma) =-\phi(\sigma, \gamma(\sigma))\) for \(\gamma\). It is left to find the initial conditions in the physical space from \eqref{CG} as
\begin{equation}
\begin{aligned}
        &\eta_0(x)=\psi(\sigma, \gamma(\sigma))-\phi^2(\sigma, \gamma(\sigma))/2,\quad u_0(x)=\phi(\sigma, \gamma(\sigma)),
        \\ 
      &x= \sigma - \psi(\sigma, \gamma(\sigma))+\phi^2(\sigma, \gamma(\sigma))/2.
\end{aligned}
\end{equation}
The uniqueness of this process follows from the invertibility of the CGT.

In conclusion, the the inverse problem can be solved exactly without using the method of data projection, which as can be seen in \eqref{eq:proj} is quite involved. Moreover, we are able to recover both the initial displacement $\eta_0$ and initial velocity $u_0$ from just the shoreline oscillations $R(t)$, which is counter-intuitive at first sight.

\section{Remarks}\label{sec:remarks}
In this section we point out several important assumptions and details. In our derivations we assume the CGT to be invertible. Mathematically it means that the Jacobians
\begin{equation}
    \det \frac{\partial (\sigma,\tau)}{\partial(x,t)},\quad 
    \det \frac{\partial(x,t)}{\partial(\sigma,\tau)} 
\end{equation}
do not vanish. Physically it means that the wave does not break \cite{Rybkin14}.

We assume \(R (t)\) to be an analytic function. In practice, the shoreline data would be given as a set of measurements. In that case analytic approximation must be employed first. For example, a wave-like data can be approximated as a finite sum of Gaussian pulses:
\begin{equation}
     \ol R (t) = \sum a_j e^{-b_j(t-c_j)^2},
\end{equation}
where \(a_j, c_j \in \R\) and \(b_j > 0\).

When we recover the curve \(\gamma(\sigma)\) in the hodograph plane we first recover \(\psi, \phi\) in the entire \((\sigma, \tau)\) plane. While this makes perfect sense mathematically, it is impossible to do numerically. Instead one should estimate how far \(\gamma(\sigma)\) is from \(\tau = 0\) based on the computed \(\phi_{\text{proj}}\), next compute \(\psi, \phi\) on a grid \((\sigma_j, \tau_i)\), and then for each \(\sigma_j\) find \(\tau_{i_j}\) that minimises \(f(\tau)=\abs{ \tau + \phi(\sigma_j,\tau)}\). Then set \(\gamma(\sigma_j):= \tau_{i_j}\). Higher accuracy can be achieved by refining the grid around the curve found and repeating this process.

Finally, we note that the integral transforms derived in \eqref{finalPsi} and \eqref{finalPhi} are a special form of the Abel transform, which finds its applications in various fields, such as image processing, tomography, and astronomy, especially in context of inverse problems \cite{Epstein07,Deans,deHoop2017}. In our derivations it appears as the composition of the Fourier and Hankel transform; this result is known as the projection-slice theorem \cite{Bracewell2003}.  
\section{Verification}\label{sec:num}
\begin{figure}[h!]
    \begin{subfigure}{.45\linewidth}
        \includegraphics[width = \linewidth]{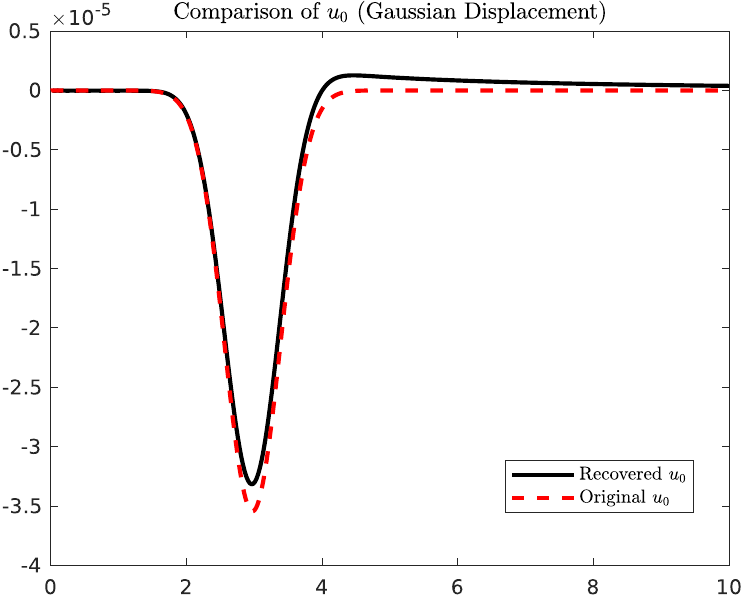}
    \end{subfigure}\hfill
    \begin{subfigure}{.45\linewidth}
        \includegraphics[width = \linewidth]{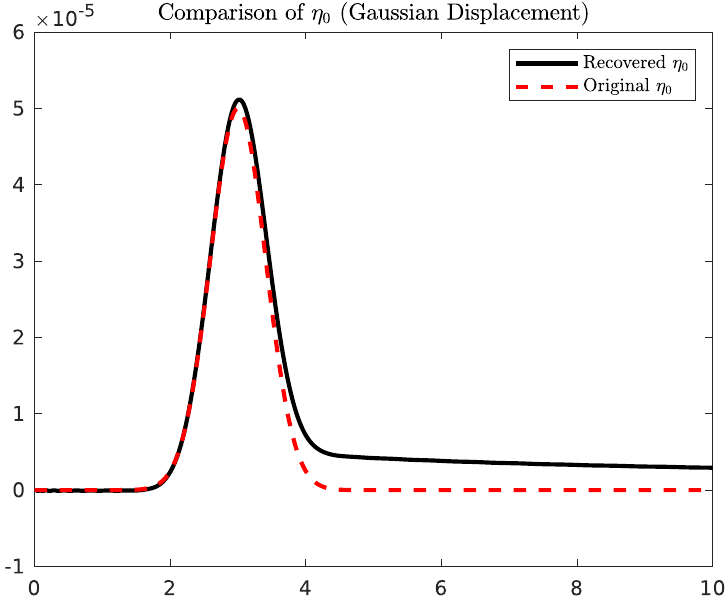}
    \end{subfigure}
    \caption{Comparison of the original initial conditions to the recovered initial conditions for a Gaussian wave.}\label{fig:icsGauss}
\end{figure}

\begin{figure}[h!]
    \begin{subfigure}{.45\linewidth}
        \includegraphics[width = \linewidth]{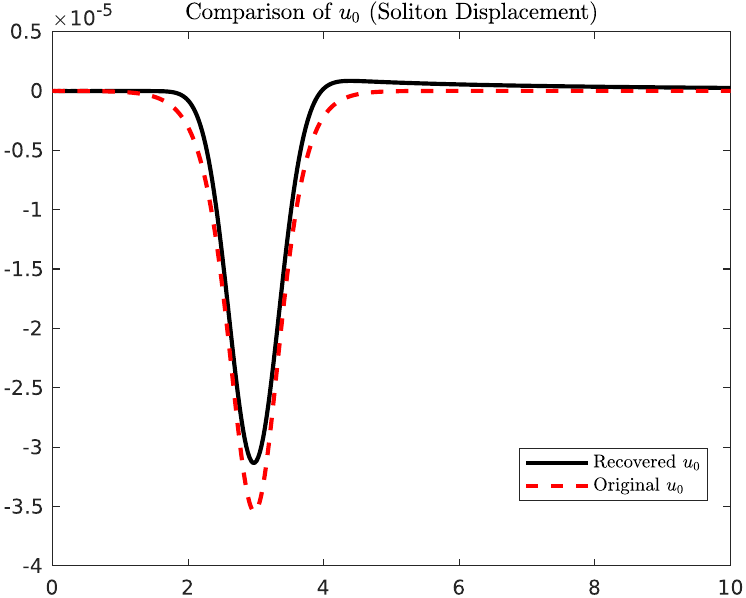}
    \end{subfigure}\hfill
    \begin{subfigure}{.45\linewidth}
        \includegraphics[width = \linewidth]{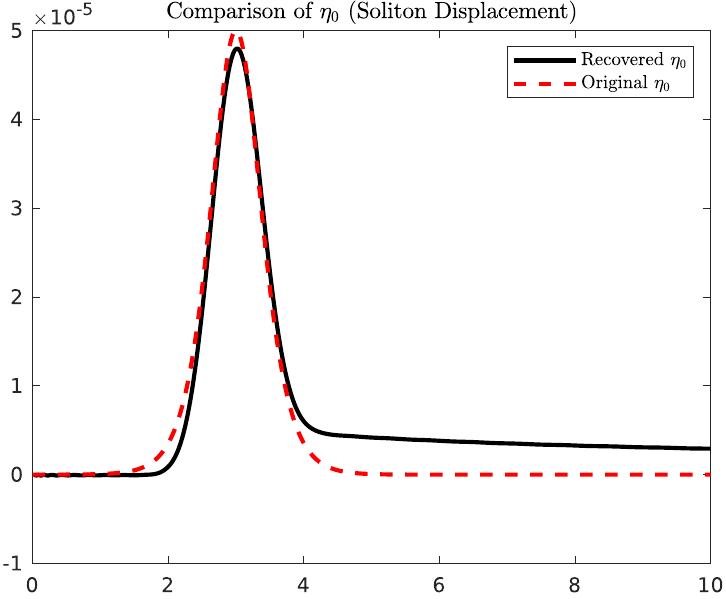}
    \end{subfigure}
    \caption{Comparison of the original initial conditions to the recovered initial conditions for a soliton wave.}\label{fig:icsS}
\end{figure}

\begin{figure}[h!]
    \begin{subfigure}{.45\linewidth}
        \includegraphics[width = \linewidth]{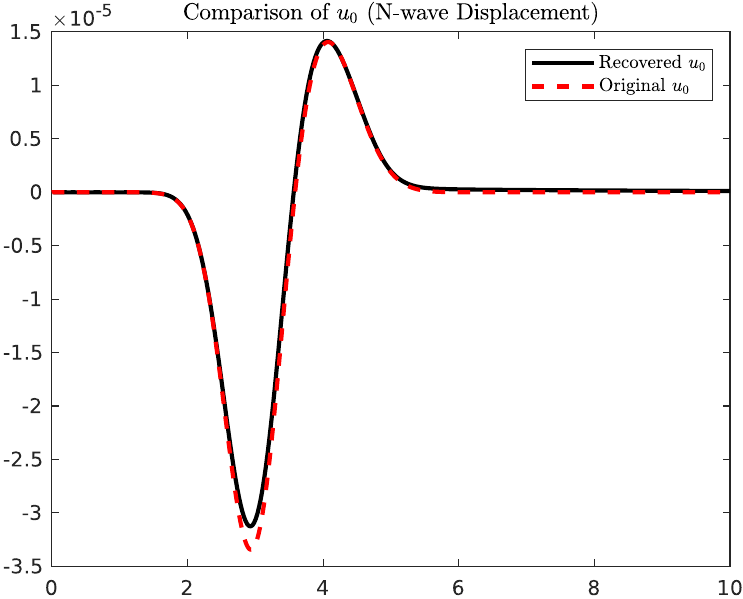}
    \end{subfigure}\hfill
    \begin{subfigure}{.45\linewidth}
        \includegraphics[width = \linewidth]{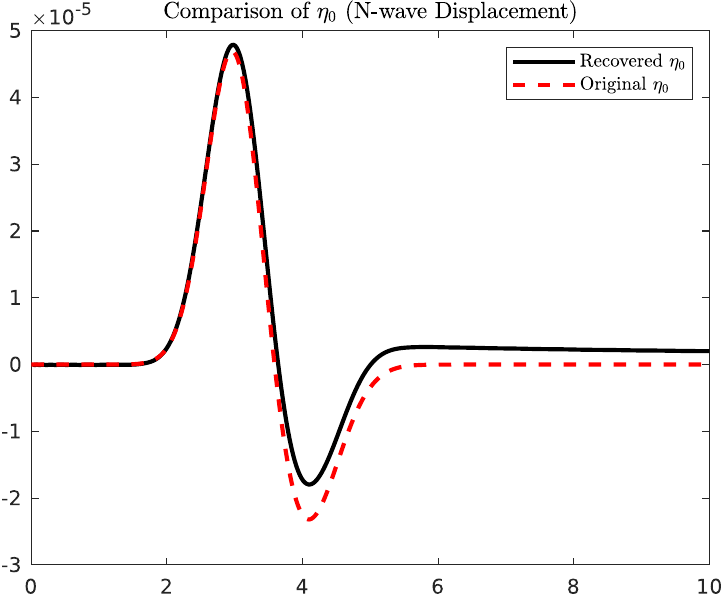}
    \end{subfigure}
    \caption{Comparison of the original initial conditions to the recovered initial conditions for an N-wave.}\label{fig:icsN}
\end{figure}

For numerical verification we manufacture the solution to the inverse problem. We take \(R (t)\) generated via the forward solution from three different initial conditions, a Gaussian pulse, a soliton wave,  and an N-wave,
\begin{equation}\label{eq:initialdisplacement}
\begin{aligned}
    &\eta_0 (x) = 0.00005e^{-3(x - 3)^2}, \\
    &\eta_0 (x) = 0.00005 \operatorname{sech} ^2(2x-6),\\
    &\eta_0 (x) = 0.00005e^{-3(x - 3)^2} - 0.000025e^{-2(x - 4)^2},
\end{aligned}
\end{equation}
and in all three cases we assume the velocity to be
\begin{equation}
    u_0 (x) = -2\sqrt{\frac{m+1}{m}}\left(\sqrt{\eta_0(x)+x} - \sqrt{x}\right) \label{velocityrelationship}.
\end{equation}
For such a velocity it was shown by \cite{Didenkulova11Velocity} that the wave propagates towards the shore.
Starting with \eqref{eq:initialdisplacement} and \eqref{velocityrelationship} the shoreline oscillations can be computed using the solution to the forward problem.
This solution to the forward problem was numerically verified 
against the so-called augmented shallow water equations system by \cite{ValianiCaleffi}. After the run-up is found from the solution to the forward problem, we find \(\psi(0,\tau)\) and \(\phi(0,\tau)\) from \eqref{CGarshore1}. The initial conditions in the hodograph plane can then be computed from \eqref{finalPsi} and \eqref{finalPhi}. After that we solve the forward problem by \eqref{psisol} and \eqref{phisol}. Finally, following the method described in Section \ref{sec:remarks} we recover the curve \(\gamma(\sigma)\) and the initial conditions on it. Figures \ref{fig:icsGauss}, \ref{fig:icsS}, and \ref{fig:icsN} demonstrate the original and the recovered initial conditions. We see that our model captures the key features of the wave.
We observe that our algorithm performs better on the N-wave, than on the Gaussian or soliton wave. We believe, this happens due to the numerical error build-up, which is more noticeable in cases of sign-constant waves (Gaussian and soliton).

While our results here are mainly analytical, we briefly discuss the computational cost of the proposed algorithm.

Given the run-up data \(R(t)\) as an \(N\)-dimensional vector of measurements, we find \(\psi(0, \tau)\) and \(\phi(0,\tau)\) from \eqref{CGarshore1} using numerical differentiation and pointwise algebraic operations with arrays (vectors).
For our next step we implement numerical integration for equations \eqref{finalPsi} and \eqref{finalPhi} in order to find \(\psi(\sigma, 0)\) and \(\phi(\sigma, 0)\). While the current quadrature is \(O(N^2)\), these integral transforms are compositions of the Fourier and Hankel transforms, so the integration can be implemented using the FFT and FHT at the cost \(O(N \log N)\).
Then, in order to perform the inverse CGT we need to compute \(\psi\) and \(\phi\) from \eqref{psisol} and \eqref{phisol}. This again can be done using the FFT and FHT for all \(\sigma, \tau\), so the computational cost would be \(O(N^2 \log N)\), but the current slow implementation is \(O(N^4)\). Next, we recover the curve \(\gamma\). Using brute force to solve the minimisation problem for \(\gamma\) described in Section \ref{sec:remarks} is \(O(N^2)\). Finally, once \(\gamma\) is found, the inverse CGT is performed using pointwise algebraic operations with vectors. Thus, if FFT and FHT algorithms are applied then the total cost is \(O(N^2\log N)\).

\section{Conclusions}\label{sec:conc}
We have shown that for the 1D non-linear shallow water equations in power-shaped bays, there exists a relationship between the hodograph functions at the shore, which allows for the recovery of both the initial displacement and velocity from the shoreline oscillations. Moreover, in the hodograph plane the initial displacement and velocity give distinguishable contributions to the run-up.
This result is surprising, as it is unexpected that both the intial conditions can be recovered from a single function \(R(t)\). 
Numerical experiments demonstrate good agreement between the original and the recovered initial condition.  Finally, we have shown that the inverse problem can be solved exactly, without loss of accuracy due to the projections, unlike the direct problem.


\section*{Acknowledgements}
Alexei Rybkin, Oleksandr Bobrovnikov, Noah Palmer, and Daniel Abramowicz acknowledge  support from NSF grant DMS-2307774.
Efim Pelinovsky's research is supported by the RSF grant 24-47-02007 (section 3).
\section*{Data availability}
No data was used for the research described in the article.
\bibliographystyle{apalike}
\bibliography{text.bib}

\end{document}